\documentclass[12pt, reqno]{amsart}
\usepackage{amsmath, amsthm, amscd, amsfonts, amssymb, graphicx, color}
\usepackage[bookmarksnumbered, colorlinks, plainpages]{hyperref}

\newtheorem{theorem}{Theorem}[section]
\newtheorem{lemma}[theorem]{Lemma}
\newtheorem{proposition}[theorem]{Proposition}
\newtheorem{corollary}[theorem]{Corollary}
\theoremstyle{definition}
\newtheorem{definition}[theorem]{Definition}

\theoremstyle{remark}
\newtheorem{remark}[theorem]{Remark}
\numberwithin{equation}{section}

\allowdisplaybreaks
\begin{document}

\title[On perturbation of continuous frames in Hilbert $C^{\ast}$-modules]
{On perturbation of continuous frames in Hilbert $C^{\ast}$-modules}

\author[H.Ghasemi]{Hadi Ghasemi}
\address{Hadi Ghasemi \\ Department of Mathematics and Computer
Sciences, Hakim Sabzevari University, Sabzevar, P.O. Box 397, IRAN}
\email{ \rm h.ghasemi@hsu.ac.ir}
\author[T.L. Shateri]{Tayebe Lal Shateri }
\address{Tayebe Lal Shateri \\ Department of Mathematics and Computer
Sciences, Hakim Sabzevari University, Sabzevar, P.O. Box 397, IRAN}
\email{ \rm  t.shateri@hsu.ac.ir; shateri@ualberta.ca}
\thanks{*The corresponding author:
t.shateri@hsu.ac.ir; shateri@ualberta.ca (Tayebe Lal Shateri)}
 \maketitle

\begin{abstract}
In the present paper, we examine the perturbation of continuous frames and Riesz-type frames in Hilbert $C^*$-modules. We extend the Casazza-Christensen general perturbation theorem for Hilbert space frames to continuous frames in Hilbert $C^*$-modules. We obtain a necessary condition under which the perturbation of a Riesz-type frame of Hilbert $C^*$-modules remains to be a Riesz-type frame. Also, we examine the effect of duality on the perturbation of continuous frames in Hilbert $C^*$-modules, and we prove if the operator frame of a continuous frame $F$ is near to the combination of the synthesis operator of a continuous Bessel mapping $G$ and the analysis operator of $F$, then $G$ is a continuous frame.
\vskip 3mm
\end{abstract}

\section{Introduction And Preliminaries}
Frame theory is nowadays a fundamental research area in mathematics,
computer science and engineering with many exciting and interesting  applications in a variety of different fields. Frames for Hilbert spaces have natural generalizations in Hilbert $C^*$-modules that are generalizations of Hilbert spaces by allowing the inner product to take values in a more general $C^*$-algebra than $\mathbb C$. 
One reason to study frames in Hilbert $C^*$-modules is that there are some differences between Hilbert spaces and Hilbert $C^*$-modules. For example, in general, every bounded operator on a Hilbert space has an unique adjoint, while this fact not hold for bounded operators on a Hilbert $C^*$-module.\\
Perturbation questions for frames in Hilbert spaces have been studied intensively in the last decade, with satisfying results for general frames \cite{CHR1}. Some results have also been obtained for frame sequences, see \cite{CHR2,CHR3}. In \cite{HA}, the authors investigate  the perturbation of frames and Riesz bases in Hilbert $C^*$-modules.

In this paper, we extend the Casazza-Christensen general perturbation theorem for Hilbert space frames to continuous frames in Hilbert $C^*$-modules. In the Hilbert space setting, under the same perturbation condition, the perturbation of any Riesz basis remains to be a Riesz basis. However, this result is no longer true for Riesz bases in Hilbert 
$C^*$-modules. We obtain a necessary condition under which the perturbation of a Riesz-type frame of Hilbert $C^*$-modules remains to be a Riesz-type frame. Also, we show that if the difference between a weakly measurable mapping $G$ with a continuous frame is a continuous Bessel mapping, then $G$ is a continuous frame. Finally, we examine the effect of duality on the perturbation of continuous frames in Hilbert $C^*$-modules, and we prove if the operator frame of a continuous frame $F$ is near to the combination of the synthesis operator of a continuous Bessel mapping $G$ and the analysis operator of $F$, then $G$ is a continuous frame.

In the following we briefly recall the definitions and basic properties of Hilbert $C^*$-modules. we refer to  \cite{LAN,OLS}. Throughout this paper, $\mathcal A$ shows a unital $C^*$-algebra.
%-------------------------------------------------------------------------------------------------------------------------------------------------------
\begin{definition}
Let $U$ be a left $\mathcal A$-module equipped with an $\mathcal A$-valued inner product $\langle .,.\rangle :U\times U\to \mathcal A$ which is $\mathbb C$-linear and $\mathcal A$-linear in its first variable and satisfies the following conditions:\\
$(i)\; \langle f,f\rangle \geq 0$,\\
$(ii)\; \langle f,f\rangle =0$  iff $f=0$,\\
$(iii)\; \langle f,g\rangle ^*=\langle g,f\rangle ,$\\
$(iv)\; \langle af,g\rangle=a\langle f,g\rangle ,$\\
for all $f,g\in U$ and $a\in\mathcal A$.
Then the pair $(H,\langle .,.\rangle)$ is called a \textit{pre-Hilbert module} over $C^*$-algebra $\mathcal A$. If $U$ is complete with respect to the topology determined by the norm $\|f\|=\|\langle f,f\rangle \|^{\frac{1}{2}}$, then it is called a \textit{Hilbert $\mathcal A$-module}.
\end{definition}
%---------------------------------------------------------------------------------------------------------------------------------------------------------
The \textit{Cauchy-Schwartz inequality} reconstructed in Hilbert $C^*$-modules as follow \cite{OLS}.
%---------------------------------------------------------------------------------------------------------------------------------------------------------
\begin{lemma}
Let $U$ be a Hilbert $C^*$-modules over a unital $C^*$-algebra $\mathcal A$. Then
\begin{equation*}
\|\langle f,g\rangle \|^{2}\leq\|\langle f,f\rangle \|\;\|\langle g,g\rangle \|
\end{equation*}
for all $f,g\in U$.
\end{lemma}
%---------------------------------------------------------------------------------------------------------------------------------------------------------
\begin{proposition}\label{NO}
\cite[Proposition 1.2.4]{MT} Let $U$ be a Hilbert $C^*$-modules over a unital $C^*$-algebra $\mathcal A$. Then\\
(i)\;$\| af\|\leq\| a\|\| f\|$\\
(ii)\;$\langle f,g\rangle\langle g,f\rangle\leq\| g\|^{2}\langle f,f\rangle$\\
for all $f,g\in U$ and $a\in\mathcal A$.
\end{proposition}
%---------------------------------------------------------------------------------------------------------------------------------------------------------
For two Hilbert $C^*$-modules $U$ and $V$ over a unital $C^*$-algebra $\mathcal A$. A map $T:U\to V$ is said to be \textit{adjointable} if there exists a map $T^{*}:V\to U$ satisfying
$$\langle Tf,g\rangle =\langle f,T^*g\rangle $$
for all $f\in U, g\in V$. Such a map $T^*$ is called the \textit{adjoint} of $T$. By $End_{\mathcal A}^*(U)$ we denote the set of all adjointable maps on $U$. It is surprising that every adjointable operator is automatically linear and bounded.
%---------------------------------------------------------------------------------------------------------------------------------------------------------
\begin{theorem}\label{BL}
\cite[Theorem 2.1.4]{MT} Let $U$ and $V$ be two Hilbert $C^*$-modules over a unital $C^*$-algebra $\mathcal A$ and $T\in End_{\mathcal A}^*(U,V)$. Then The following are equivalent:\\
(i)\;$T$ is bounded and $\mathcal A$-linear.\\
(ii)\;There exists $k>0$ such that $\langle Tf,Tf\rangle\leq k\langle f,f\rangle$, for all $f\in U$.
\end{theorem}
%-------------------------------------------------------------------------------------------------------------------------------------------
\begin{lemma}\label{SB} 
\cite{AR} Let $U$ and $V$ be two Hilbert $C^*$-modules over a unital $C^*$-algebra $\mathcal A$ and $T\in End_{\mathcal A}^*(U,V)$. Then The following are equivalent:\\
$(i)$\;$T$ is surjective\\
$(ii)$\;$T^{*}$ is bounded below with respect to the norm, i.e. there exists $m>0$ such that\\ $\|T^{*}f\|\geq m\| f\|$, for all $f\in V$.\\
$(iii)$\;$T^{*}$ is bounded below with respect to the inner product, i.e. there exists $m>0$ such that $\langle T^{*}f,T^{*}f\rangle\geq m\langle f,f\rangle$, for all $f\in V$.
\end{lemma}
%-------------------------------------------------------------------------------------------------------------------------------------------
Let $\mathcal Y$ be a Banach space, $(\mathcal X,\mu)$ a measure space, and $f:\mathcal X\to \mathcal Y$ a measurable
function. The integral of the Banach-valued function $f$ has been defined by Bochner and others. Most properties of this integral are similar to those of the integral of real-valued functions (see \cite{DAN,YOS}). Since every $C^*$-algebra and Hilbert $C^*$-module is a Banach space, hence we can use this integral in these spaces.
%-------------------------------------------------------------------------------------------------------------------------------------------
In the sequel, we assume that $\mathcal A$ is a unital $C^*$-algebra and $U$ is a Hilbert $C^*$-module over $\mathcal A$ and $(\Omega ,\mu)$ is a measure space.
%-------------------------------------------------------------------------------------------------------------------------------------------
Define,
\begin{equation*}
L^{2}(\Omega ,A)=\lbrace\varphi :\Omega \to A\quad ;\quad \Vert\int_{\Omega}\vert(\varphi(\omega))^{*}\vert^{2} d\mu(\omega)\Vert<\infty\rbrace.
\end{equation*}
It was shown that $L^{2}(\Omega , A)$ is a Hilbert $\mathcal A$-module with the inner product $$\langle \varphi ,\psi\rangle = \int_{\Omega}\langle\varphi(\omega),\psi(\omega)\rangle d\mu(\omega),\quad (\varphi ,\psi \in L^{2}(\Omega , A)),$$ and induced norm $\|\varphi\|=\|\langle \varphi,\varphi\rangle \|^{\frac{1}{2}}$, \cite{LAN}.\\
In the following, we recall the notion of continuous frames in Hilbert $C^*$-modules over a unital $C^*$-algebra $\mathcal A$, and some properties of these frames are given. For details, see \cite{GHSH1,GHSH}.
%-------------------------------------------------------------------------------------------------------------------------------------------
\begin{definition}
A mapping $F:\Omega \to U$ is called a continuous frame for $U$ if\\
$(i)\; F$ is weakly-measurable, i.e, the mapping $\omega\longmapsto\langle f,F(\omega)\rangle $ is measurable on $\Omega$, for any $f\in U$.\\
$(ii)$ There exist constants $A,B>0$ such that
\begin{equation}\label{eq1}
A\langle f,f\rangle \leq \int_{\Omega}\langle f,F(\omega)\rangle \langle F(\omega),f\rangle d\mu(\omega)\leq B\langle f,f\rangle  ,\quad (f\in U). 
\end{equation}
\end{definition}
%-------------------------------------------------------------------------------------------------------------------------------------------
The constants $A,B$ are called the \textit{lower} and the \textit{upper} frame bounds, respectively. The mapping $F$ is called \textit{Bessel} if the right inequality in \eqref{eq1} holds and is called \textit{tight} if  $A=B$.\\
%-------------------------------------------------------------------------------------------------------------------------------------------
For a continuous frame $F:\Omega \to U$ the following operators were defined.\\
$(i)$\;The \textit{synthesis operator} or \textit{pre-frame operator} $T_{F}:L^{2}(\Omega , A)\;\to U$ weakly defined by
\begin{equation}
\langle T_{F}\varphi ,f\rangle =\int_{\Omega}\varphi(\omega)\langle F(\omega),f\rangle d\mu(\omega),\quad (f\in U).
\end{equation}
$(ii)$\; The adjoint of $T$, called the \textit{analysis operator} $T^{\ast}_{F}:U\;\to L^{2}(\Omega , A)$ is defined by
\begin{equation}
(T^{\ast}_{F}f)(\omega)=\langle f ,F(\omega)\rangle\quad (\omega\in \Omega).
\end{equation}
$(iii)$\; The \textit{frame operator} $S_{F}:U\;\to U$ is weakly defined by
\begin{equation}
\langle S_{F}f ,f\rangle =\int_{\Omega}\langle f,F(\omega)\rangle\langle F(\omega),f\rangle d\mu(\omega),\quad (f\in U).
\end{equation}
In \cite{GHSH} was proved the pre-frame operator $T_{F}:L^{2}(\Omega , A)\;\to U$ is well defined, surjective, adjointable $\mathcal A$-linear map and bounded with $\| T\|\leq\sqrt{B}$ . Moreover the analysis operator $T^{\ast}_{F}:U\;\to L^{2}(\Omega , A)$ is injective and has closed range. Also $S=T T^{\ast}$ is positive, adjointable, self-adjoint and invertible and $\| S\|\leq B$.\\
%-------------------------------------------------------------------------------------------------------------------------------------------
Now we recall the concept of duals of continuous frames in Hilbert $C^{\ast}$-modules \cite{GHSH3}.
 %-------------------------------------------------------------------------------------------------------------------------------------------
\begin{definition}
Let $F:\Omega \to U$ be a continuous Bessel mapping. A continuous Bessel mapping $G:\Omega \to U$ is called a \textit{dual} for $F$ if
\begin{equation*}
f= \int_{\Omega}\langle f,G(\omega)\rangle F(\omega)d\mu(\omega),\qquad(f\in U)
\end{equation*}
or
\begin{equation}\label{eq2}
\langle f,g\rangle= \int_{\Omega}\langle f,G(\omega)\rangle\langle F(\omega),g\rangle d\mu(\omega),\qquad(f,g\in U).
\end{equation}
In this case $(F,G)$ is called a \textit{dual pair}. If $T_{F}$ and $T_{G}$ denote the synthesis operators of $F$ and $G$, respectively, then \eqref{eq2} is equivalent to $T_{F}T^{*}_{G}=I_{U}$. 
\end{definition}
%-------------------------------------------------------------------------------------------------------------------------------------------
\begin{remark}
Let $F:\Omega \to U$ be a continuous frame for Hilbert $C^{\ast}$-module $U$. Then by reconstructin formula we have
\begin{equation*}
f=\int_{\Omega}\langle f,F(\omega)\rangle S^{-1}F(\omega) d\mu(\omega)=\int_{\Omega}\langle f,S^{-1}F(\omega)\rangle F(\omega) d\mu(\omega).
\end{equation*}
\end{remark}
%-------------------------------------------------------------------------------------------------------------------------------------------
Therefore $S^{-1}F$ is a dual for $F$, which is called \textit{canonical dual}. If a continuous frame $F$ has only one dual, $F$ is called a \textit{Riesz-type frame}.
%-------------------------------------------------------------------------------------------------------------------------------------------
Also we recall the concept of continuous Riesz basis in Hilbert $C^{\ast}$-modules \cite{GHSH2}.
%-------------------------------------------------------------------------------------------------------------------------------------------
\begin{definition}
Let $U$ be a Hilbert $C^{\ast}$-module $U$ over a unital $C^*$-algebra $\mathcal A$. A weakly-measurable mapping $F:\Omega \to U$ is called a \textit{continuous Riesz basis} for Hilbert $C^{\ast}$-module $U$, if the following conditions are satisfied:\\
$(i)\;$ $F$ is $\mu$-complete.\\
$(ii)\;$ There are two constants $A,B>0$ such that
\begin{equation}
A\|\int_{\Omega_{1}}\vert \varphi(\omega)^{*}\vert^{2} d\mu(\omega)\|^{\dfrac{1}{2}} \leq \|\int_{\Omega_{1}}\varphi(\omega)F(\omega) d\mu(\omega)\|\leq B\|\int_{\Omega_{1}}\vert \varphi(\omega)^{*}\vert^{2}d\mu(\omega)\|^{\dfrac{1}{2}},
\end{equation}
for every $\varphi\in L^{2}(\Omega , A)$ and measurable subset $\Omega_{1}\subseteq\Omega$ with  $\mu(\Omega_{1})<+\infty$.
\end{definition}
%-------------------------------------------------------------------------------------------------------------------------------------------
We need the following theorems in the next section.
%-------------------------------------------------------------------------------------------------------------------------------------------
\begin{theorem}\label{c-f-norm}
\cite[Theorem 2.14]{GHSH} Let $U$ be a Hilbert $C^{\ast}$-module over a unital $C^*$-algebra $\mathcal A$. Then the following are equivalent:\\
(i)\;$F:\Omega \to U$ is a continuous frame for $U$.\\
(ii)\;The mapping $\Omega \to \langle f,F(\Omega)\rangle$ is measurable and there exsist constants $A,B>0$ such that
\begin{equation*}
A\Vert \langle f,f\rangle\Vert\leq\Vert\int_{\Omega}\langle f,F(\omega)\rangle\langle F(\omega),f\rangle d\mu(\omega)\Vert\leq B\Vert \langle f,f\rangle\Vert,
\end{equation*}
for all $f\in U$
\end{theorem}
%-------------------------------------------------------------------------------------------------------------------------------------------
\begin{theorem}\label{RBT}
\cite{GHSH2}Let $F:\Omega \to U$ be a continuous frame for Hilbert $C^{\ast}$-module $U$ over a unital $C^*$-algebra $\mathcal A$. Then $F$ is a continuous Riesz basis for $U$ if and only if $F$ is a Riesz-type frame.
\end{theorem}
%-------------------------------------------------------------------------------------------------------------------------------------------
\begin{theorem}
\cite[Theorem 3.4]{GHSH} Let $F:\Omega \to U$ be a continuous frame for Hilbert $C^{\ast}$-module $U$ over a unital $C^*$-algebra $\mathcal A$. Then $F$ is a Riesz-type frame if and only if the analysis operator $T^{*}_{F}:U\to L^{2}(\Omega ,A)$ is onto.
\end{theorem}
%-------------------------------------------------------------------------------------------------------------------------------------------
\begin{lemma}\label{aF-B}
\cite{GHSH4}Let $F:\Omega \to U$ be a continuous Bessel mapping for Hilbert $C^{\ast}$-module $U$ over a unital $C^*$-algebra $\mathcal A$ with the bound $B$ and $a\in \mathcal A$. Then $aF:\Omega \to U$ is a continuous Bessel mapping for $U$.
\end{lemma}
%------------------------------------------------------------------------------------------------------------------------------------------- 
\begin{remark}\label{aF-f}
\cite{GHSH4}Let $F:\Omega \to U$ be a continuous frame for Hilbert $C^{\ast}$-module $U$ over a unital $C^*$-algebra $\mathcal A$ and $a\in \mathcal A$. It is easily to see that:\\
$(1)$\; If $a$ is unitary, then $aF$ is a continuous frame for $U$ and $S_{aF}=S_{F}$.\\
$(2)$\; If $a$ is unitary and $G:\Omega \to U$ is a dual of $F$, then $aG$ is a dual of $aF$.\\
$(3)$\; If  $ab=ba$, for all $b\in \mathcal A$, then $aF$ is a continuous frame for $U$ and $S_{aF}=\vert a\vert^{2}S_{F}$ and $T_{aF}=aT_{F}$.
\end{remark}
%-------------------------------------------------------------------------------------------------------------------------------------------
%------------------------------------------------------------------------------------------------------------------------------------------- 
%-------------------------------------------------------------------------------------------------------------------------------------------
%------------------------------------------------------------------------------------------------------------------------------------------- 
%-------------------------------------------------------------------------------------------------------------------------------------------
%------------------------------------------------------------------------------------------------------------------------------------------- 
\section{Main results}

In this section, we examine the perturbation of continuous frames and Riesz-type frames in Hilbert $C^*$-modules.

First, we show the condition that the sum of a continuous frame and a Bessel mapping is a continuous frame.
%-------------------------------------------------------------------------------------------------------------------------------------------
\begin{theorem}\label{sum3}
Let $F:\Omega \to U$ be a continuous frame for Hilbert $C^{\ast}$-module $U$ over a unital $C^*$-algebra $\mathcal A$ with the frame bounds $A,B>0$. Also assume that $a_{1},a_{2}\in \mathcal A$ are invertible and $a_{i}b=ba_{i}$, for $i=1,2$ and all $b\in \mathcal A$. If $G:\Omega \to U$ is a continuous Bessel mapping for $U$ with bound $N$ such that $N\Vert a_{2}\Vert^{2}<A\Vert a_{1}^{-1}\Vert^{-2}$, then $a_{1}F+a_{2}G$ is a continuous frame for $U$.
\end{theorem}

\begin{proof}
For any $f\in U$, we obtain
\begin{align*}
\Vert\int_{\Omega}\langle f,& a_{1}F(\omega)+a_{2}G(\omega)\rangle\langle a_{1}F(\omega)+a_{2}G(\omega),f\rangle d\mu(\omega)\Vert^{\frac{1}{2}}=\Vert\lbrace\langle f,a_{1}F(\omega)+a_{2}G(\omega)\rangle\rbrace_{\omega\in\Omega}\Vert\\
& \leq\Vert\lbrace\langle f,a_{1}F(\omega)\rangle\rbrace_{\omega\in\Omega}\Vert +\Vert\lbrace\langle f,a_{2}G(\omega)\rangle\rbrace_{\omega\in\Omega}\Vert\\
& =\Vert\int_{\Omega}\langle f,F(\omega)\rangle a_{1}^{*}a_{1}\langle F(\omega),f\rangle d\mu(\omega)\Vert^{\frac{1}{2}}+\Vert\int_{\Omega}\langle f,G(\omega)\rangle a_{2}^{*}a_{2}\langle G(\omega),f\rangle d\mu(\omega)\Vert^{\frac{1}{2}}\\
& \leq\Vert a_{1}\Vert\Vert\int_{\Omega}\langle f,F(\omega)\rangle\langle F(\omega),f\rangle d\mu(\omega)\Vert^{\frac{1}{2}}+\Vert a_{2}\Vert\Vert\int_{\Omega}\langle f,G(\omega)\rangle\langle G(\omega),f\rangle d\mu(\omega)\Vert^{\frac{1}{2}}\\
& \leq\Vert a_{1}\Vert\sqrt{B}\Vert\langle f,f\rangle\Vert^{\frac{1}{2}}+\Vert a_{2}\Vert\sqrt{N}\Vert\langle f,f\rangle\Vert^{\frac{1}{2}}\\
& =(\Vert a_{1}\Vert\sqrt{B}+\Vert a_{2}\Vert\sqrt{N})\Vert\langle f,f\rangle\Vert^{\frac{1}{2}}.
\end{align*}
Also,
\begin{align*}
\Vert\int_{\Omega}\langle f,& a_{1}F(\omega)+a_{2}G(\omega)\rangle\langle a_{1}F(\omega)+a_{2}G(\omega),f\rangle d\mu(\omega)\Vert^{\frac{1}{2}}=\Vert\lbrace\langle f,a_{1}F(\omega)+a_{2}G(\omega)\rangle\rbrace_{\omega\in\Omega}\Vert\\
& \geq\Vert\lbrace\langle f,a_{1}F(\omega)\rangle\rbrace_{\omega\in\Omega}\Vert -\Vert\lbrace\langle f,a_{2}G(\omega)\rangle\rbrace_{\omega\in\Omega}\Vert\\
& =\Vert\int_{\Omega}\langle f,F(\omega)\rangle a_{1}^{*}a_{1}\langle F(\omega),f\rangle d\mu(\omega)\Vert^{\frac{1}{2}}-\Vert\int_{\Omega}\langle f,G(\omega)\rangle a_{2}^{*}a_{2}\langle G(\omega),f\rangle d\mu(\omega)\Vert^{\frac{1}{2}}\\
& \geq\Vert a_{1}^{-1}\Vert^{-1}\sqrt{A}\Vert\langle f,f\rangle\Vert^{\frac{1}{2}}-\Vert a_{2}\Vert\sqrt{N}\Vert\langle f,f\rangle\Vert^{\frac{1}{2}}\\
& =(\Vert a_{1}^{-1}\Vert^{-1}\sqrt{A}-\Vert a_{2}\Vert\sqrt{N})\Vert\langle f,f\rangle\Vert^{\frac{1}{2}}.
\end{align*}
Hence $a_{1}F+a_{2}G$ is a continuous frame for $U$, by Theorem \ref{c-f-norm}.
\end{proof}
%-------------------------------------------------------------------------------------------------------------------------------------------
Replacing $G$ by $a_{2}^{-1}G-a_{2}^{-1}a_{1}F$ in Theorem \ref{sum3}, we obtain the next corollary.
%-------------------------------------------------------------------------------------------------------------------------------------------
\begin{corollary}\label{sum4}
Let $F:\Omega \to U$ be a continuous frame for Hilbert $C^{\ast}$-module $U$ over a unital $C^*$-algebra $\mathcal A$ with the frame bounds $A,B>0$. Also assume that $a_{1},a_{2}\in \mathcal A$ are invertible and $a_{i}b=ba_{i}$, for $i=1,2$ and all $b\in \mathcal A$. If $G:\Omega \to U$ is a weakly measurable mapping such that $a_{2}^{-1}G-a_{2}^{-1}a_{1}F$ is a continuous Bessel mapping for $U$ with bound $N$ and $N\Vert a_{2}\Vert^{2}<A\Vert a_{1}^{-1}\Vert^{-2}$, then $G$ is a continuous frame for $U$.
\end{corollary}
%-------------------------------------------------------------------------------------------------------------------------------------------
Setting $a_{1}=-a_{2}=1_{\mathcal A}$, the identity element of $\mathcal A$,we have the next corollary.
%-------------------------------------------------------------------------------------------------------------------------------------------
\begin{corollary}\label{pert(F-G)B}
Let $F:\Omega \to U$ be a continuous frame for Hilbert $C^{\ast}$-module $U$ over a unital $C^*$-algebra $\mathcal A$ with the lower bound $A$. If $G:\Omega \to U$ is a weakly measurable mapping such that $F-G$ is a continuous Bessel mapping for $U$ with bound $N$ and $N<A$, then $G$ is a continuous frame for $U$.
\end{corollary}
%-------------------------------------------------------------------------------------------------------------------------------------------
In the following we investigate some conditions under which a mapping becomes a continuous frame under the influence of a continuous frame.
%-------------------------------------------------------------------------------------------------------------------------------------------
\begin{theorem}\label{pert1}
Let $F:\Omega \to U$ be a continuous frame for Hilbert $C^{\ast}$-module $U$ over a unital $C^*$-algebra $\mathcal A$ with the frame bounds $A,B>0$. Also assume that $a_{1},a_{2}\in \mathcal A$ are invertible and $a_{i}b=ba_{i}$, for $i=1,2$ and all $b\in \mathcal A$. If $G:\Omega \to U$ is a weakly measurable mapping and there exist constants $\alpha ,\beta, \gamma\geq 0$ such that $max\lbrace\beta\;,\;\beta\Vert a_{1}^{-1}\Vert\Vert a_{1}\Vert\;,\;\Vert a_{1}^{-1}\Vert(\alpha\Vert a_{1}\Vert+\dfrac{\gamma}{\sqrt{A}})\rbrace<1$ and
\begin{align*}
\Vert\int_{\Omega}\psi(\omega)\langle a_{1}F(\omega)-a_{2}G(\omega),f\rangle d\mu(\omega)\Vert\leq & \alpha\Vert\int_{\Omega}\psi(\omega)\langle a_{1}F(\omega),f\rangle d\mu(\omega)\Vert\\
& +\beta\Vert\int_{\Omega}\psi(\omega)\langle a_{2}G(\omega),f\rangle d\mu(\omega)\Vert +\gamma\Vert\psi\Vert ,
\end{align*}
for all $\psi\in L^{2}(\Omega ,\mathcal A)$ and $f\in U$, then $G$ is a continuous frame for $U$.
\end{theorem}

\begin{proof}
Suppose that $T_{F}$ and $S_{F}$ are the pre-frame operator and the frame operator of $F$, respectively. For any $f\in U$ and $\psi\in L^{2}(\Omega ,\mathcal A)$ we have,
\begin{align*}
\Vert\int_{\Omega}\psi(\omega)\langle a_{2}G(\omega),f\rangle d\mu(\omega)\Vert &\leq\Vert\int_{\Omega}\psi(\omega)\langle a_{1}F(\omega),f\rangle d\mu(\omega)\Vert\\
&\;\;\;\; +\Vert\int_{\Omega}\psi(\omega)\langle a_{1}F(\omega)-a_{2}G(\omega),f\rangle d\mu(\omega)\Vert\\
&\leq (1+\alpha)\Vert\int_{\Omega}\psi(\omega)\langle a_{1}F(\omega),f\rangle d\mu(\omega)\Vert\\
&\;\;\;\; +\beta\Vert\int_{\Omega}\psi(\omega)\langle a_{2}G(\omega),f\rangle d\mu(\omega)\Vert +\gamma\Vert\psi\Vert .
\end{align*}
Since
\begin{align*}
\Vert\int_{\Omega}\psi(\omega)\langle a_{1}F(\omega),f\rangle d\mu(\omega)\Vert & =\Vert a_{1}\int_{\Omega}\psi(\omega)\langle F(\omega),f\rangle d\mu(\omega)\Vert\\
& \leq\Vert a_{1}\Vert\Vert\int_{\Omega}\psi(\omega)\langle F(\omega),f\rangle d\mu(\omega)\Vert\\
& =\Vert a_{1}\Vert\Vert\langle T_{F}\psi ,f\rangle\Vert\leq\Vert a_{1}\Vert\Vert T_{F}\Vert\Vert\psi\Vert\Vert f\Vert ,
\end{align*}
so
\begin{equation*}
\Vert\int_{\Omega}\psi(\omega)\langle a_{2}G(\omega),f\rangle d\mu(\omega)\Vert\leq\dfrac{(1+\alpha)\Vert a_{1}\Vert\sqrt{B}\Vert f\Vert +\gamma}{1-\beta}\Vert\psi\Vert.
\end{equation*}
Suppose that $T_{G}:L^{2}(\Omega ,\mathcal A)\longrightarrow U$ is weakly defined by
\begin{equation*}
\langle T_{G}\psi ,f\rangle =\int_{\Omega}\psi(\omega)\langle G(\omega),f\rangle d\mu(\omega),
\end{equation*}
for any $f\in U$ and $\psi\in L^{2}(\Omega ,\mathcal A)$. Then
\begin{align*}
\Vert T_{G}\psi\Vert &=\sup_{\Vert f\Vert\leq1}\Vert\langle T_{G}\psi ,f\rangle\Vert\\
&= \sup_{\Vert f\Vert\leq1}\Vert a_{2}^{-1}\Vert\Vert\int_{\Omega}\psi(\omega)\langle a_{2}G(\omega),f\rangle d\mu(\omega)\Vert\\
&\leq\dfrac{(1+\alpha)\sqrt{B}\Vert a_{1}\Vert\Vert a_{2}^{-1}\Vert +\gamma}{1-\beta}\Vert\psi\Vert.
\end{align*}
This shows that $T_{G}$ is well-defined and bounded. So by \cite[Corollary 2.17]{GHSH}, $G$ is a Bessel mapping for $U$. Also $a_{1}^{-1}a_{2}G$ is a Bessel mapping for $U$ with the pre-frame operator $T_{a_{1}^{-1}a_{2}G}=a_{1}^{-1}a_{2}T_{G}$, by Lemma \ref{aF-B}. Define $K=T_{a_{1}^{-1}a_{2}G}T^{*}_{F}S^{-1}_{F}$. Then for every $f,g\in U$,
\begin{equation*}
\langle Kf,g\rangle =\langle T_{a_{1}^{-1}a_{2}G}T^{*}_{F}S^{-1}_{F}f,g\rangle =a_{1}^{-1}\int_{\Omega}\langle f,S_{F}^{-1}F(\omega)\rangle\langle a_{2}G(\omega),g\rangle d\mu(\omega),
\end{equation*}
and
\begin{equation*}
\langle f,g\rangle =a_{1}^{-1}\int_{\Omega}\langle f,S_{F}^{-1}F(\omega)\rangle\langle a_{1}F(\omega),g\rangle d\mu(\omega).
\end{equation*}
Let $\psi_{f}:\Omega\longrightarrow \mathcal A$ be a mapping definsd by $\psi_{f}(\omega)=\langle f,S^{-1}F(\omega)\rangle$. Then $\psi_{f}\in L^{2}(\Omega ,\mathcal A)$ and $\Vert\psi_{f}\Vert\leq\dfrac{1}{\sqrt{A}}\Vert f\Vert$ for each $f\in U$.\\
Hence
\begin{equation*}
\langle f-Kf,g\rangle=a_{1}^{-1}(\int_{\Omega}\psi_{f}(\omega)\langle a_{1}F(\omega)-a_{2}G(\omega),g\rangle d\mu(\omega)),
\end{equation*}
and by assumption,
\begin{align*}
\Vert\langle f-Kf,g\rangle\Vert &=\Vert a_{1}^{-1}\int_{\Omega}\psi_{f}(\omega)\langle a_{1}F(\omega)-a_{2}G(\omega),g\rangle d\mu(\omega)\Vert\\
&\leq\Vert a_{1}^{-1}\Vert\Vert\int_{\Omega}\psi_{f}(\omega)\langle a_{1}F(\omega)-a_{2}G(\omega),g\rangle d\mu(\omega)\Vert\\
&\leq\Vert a_{1}^{-1}\Vert(\alpha\Vert a_{1}\langle f,g\rangle\Vert +\beta\Vert a_{1}\langle Kf,g\rangle\Vert+\gamma\Vert\psi_{f}\Vert)\\
&\leq\Vert a_{1}^{-1}\Vert(\alpha\Vert a_{1}\Vert\Vert\langle f,g\rangle\Vert +\beta\Vert a_{1}\Vert\Vert\langle Kf,g\rangle\Vert+\gamma\Vert\psi_{f}\Vert),
\end{align*}
for every $f,g\in U$ .So
\begin{align*}
\Vert f-Kf\Vert &= \sup_{\Vert g\Vert\leq 1}\Vert\langle f-Kf,g\rangle\Vert\\
&\leq\Vert a_{1}^{-1}\Vert(\alpha\Vert a_{1}\Vert\Vert f\Vert +\beta\Vert a_{1}\Vert\Vert Kf\Vert+\gamma\Vert\psi_{f}\Vert)\\
&\leq\Vert a_{1}^{-1}\Vert(\alpha\Vert a_{1}\Vert +\frac{\gamma}{\sqrt{A}})\Vert f\Vert +\beta\Vert a_{1}^{-1}\Vert\Vert a_{1}\Vert\Vert Kf\Vert.
\end{align*}
Hence by \cite[Lemma 4.1]{RND}, $K$ is invertible and
\begin{equation*}
\Vert K\Vert \leq\dfrac{1+\Vert a_{1}^{-1}\Vert(\alpha\Vert a_{1}\Vert +\frac{\gamma}{\sqrt{A}})}{1-\beta\Vert a_{1}^{-1}\Vert\Vert a_{1}\Vert}\;\;\;\;\;\; and\;\;\;\;\;\;\Vert K^{-1}\Vert \leq\dfrac{1+\beta\Vert a_{1}^{-1}\Vert\Vert a_{1}\Vert}{1-\Vert a_{1}^{-1}\Vert(\alpha\Vert a_{1}\Vert +\frac{\gamma}{\sqrt{A}})}.
\end{equation*}
Also for $f\in U$, we have
\begin{align*}
\langle f,f\rangle &= \langle KK^{-1}f,f\rangle\\
&=a_{1}^{-1}a_{2}\int_{\Omega}\langle K^{-1}f,S_{F}^{-1}F(\omega)\rangle\langle G(\omega),f\rangle d\mu(\omega)\\
&=a_{1}^{-1}a_{2}\langle\lbrace\langle K^{-1}f,S_{F}^{-1}F(\omega)\rangle\rbrace_{\omega\in\Omega} ,\lbrace\langle f,G(\omega)\rangle\rbrace_{\omega\in\Omega}\rangle,
\end{align*}
then
\begin{align*}
\Vert f\Vert^{4} &=\Vert\langle f,f\rangle\Vert^{2}\\
&\leq\Vert a_{1}^{-1}a_{2}\Vert\Vert\lbrace\langle K^{-1}f,S_{F}^{-1}F(\omega)\rangle\rbrace_{\omega\in\Omega}\Vert^{2}\Vert\lbrace\langle f,G(\omega)\rangle\rbrace_{\omega\in\Omega}\Vert^{2}\\
&=\Vert a_{1}^{-1}a_{2}\Vert\Vert\int_{\Omega}\langle K^{-1}f,S_{F}^{-1}F(\omega)\rangle\langle S_{F}^{-1}F(\omega),K^{-1}f\rangle d\mu(\omega)\Vert\Vert\int_{\Omega}\langle f,G(\omega)\rangle\langle G(\omega),f\rangle d\mu(\omega)\Vert\\
&= \Vert a_{1}^{-1}a_{2}\Vert\Vert\langle S^{-1}_{F}K^{-1}f,K^{-1}f\rangle\Vert\Vert\int_{\Omega}\langle f,G(\omega)\rangle\langle G(\omega),f\rangle d\mu(\omega)\Vert\\
&\leq\Vert a_{1}^{-1}a_{2}\Vert\Vert S_{F}^{-1}\Vert\Vert K^{-1}\Vert^{2}\Vert f\Vert^{2}\Vert\int_{\Omega}\langle f,G(\omega)\rangle\langle G(\omega),f\rangle d\mu(\omega)\Vert\\
&\leq\dfrac{\Vert a_{1}^{-1}a_{2}\Vert}{A}(\dfrac{1+\beta\Vert a_{1}^{-1}\Vert\Vert a_{1}\Vert}{1-\Vert a_{1}^{-1}\Vert(\alpha\Vert a_{1}\Vert +\frac{\gamma}{\sqrt{A}})})^{2}\Vert f\Vert^{2}\Vert\int_{\Omega}\langle f,G(\omega)\rangle\langle G(\omega),f\rangle d\mu(\omega)\Vert,
\end{align*}
and
\begin{equation*}
A\Vert a_{1}^{-1}a_{2}\Vert^{-1}(\frac{1-\Vert a_{1}^{-1}\Vert(\alpha\Vert a_{1}\Vert +\frac{\gamma}{\sqrt{A}})}{1+\beta\Vert a_{1}^{-1}\Vert\Vert a_{1}\Vert})^{2}\Vert\langle f,f\rangle\Vert\leq\Vert\int_{\Omega}\langle f,G(\omega)\rangle\langle G(\omega),f\rangle d\mu(\omega)\Vert.
\end{equation*}
Therefore $G$ is a continuous frame for $U$, by Theorem \ref{c-f-norm}.
\end{proof}
%-------------------------------------------------------------------------------------------------------------------------------------------
If $F:\Omega \to U$ is a Riesz-type frame, then by Theorem \ref{RBT}, $F$ is a continuous Riesz basis with the lower bound 
 $M$. We get the following result due to Theorem \ref{pert1}.
%-------------------------------------------------------------------------------------------------------------------------------------------
\begin{corollary}\label{pert2}
Let $F:\Omega \to U$ be a Riesz-type frame for Hilbert $C^{\ast}$-module $U$ over a unital $C^*$-algebra $\mathcal A$ with the frame bounds $A,B>0$ and the synthesis operator $T_{F}$. Also assume that $a_{1},a_{2}\in \mathcal A$ are invertible and $a_{i}b=ba_{i}$, for $i=1,2$ and all $b\in \mathcal A$. If $G:\Omega \to U$ is a weakly measurable mapping and there exist constants $\alpha ,\beta, \gamma\geq 0$ such that $max\lbrace \beta\Vert a_{2}\Vert\Vert a_{2}^{-1}\Vert , \beta\Vert a_{1}\Vert\Vert a_{1}^{-1}\Vert, \alpha \Vert a_{1}^{-1}\Vert^{-1}+\frac{\gamma}{M},\Vert a_{1}^{-1}\Vert(\alpha\Vert a_{1}\Vert+\dfrac{\gamma}{\sqrt{A}})\rbrace<1$ and
\begin{equation*}
\Vert\int_{\Omega}\psi(\omega)(a_{1}F(\omega)-a_{2}G(\omega)) d\mu(\omega)\Vert\leq\alpha\Vert\int_{\Omega}\psi(\omega)a_{1}F(\omega) d\mu(\omega)\Vert +\beta\Vert\int_{\Omega}\psi(\omega)a_{2}G(\omega) d\mu(\omega)\Vert +\gamma\Vert\psi\Vert ,
\end{equation*}
for all $\psi\in L^{2}(\Omega ,\mathcal A)$. Then $G$ is a Riesz-type frame for $U$.
\end{corollary}

\begin{proof}
By Theorem \ref{pert1}, the mapping $G$ is a continuous frame for $U$. Then $T^{*}_{G}$, the analysis operator of $G$, is injective. If $f\in U$ and $T^{*}_{G}f=0$, then $\langle f,G(\omega)\rangle =0$, for all $\omega\in\Omega$. This implies that $f=0$. Then $G$ is $\mu$-complete.

On the other hand, for each $\psi\in L^{2}(\Omega ,\mathcal A)$ we have
\begin{align*}
\Vert\int_{\Omega}\psi(\omega)G(\omega)d\mu(\omega)\Vert &\leq\Vert a_{2}^{-1}\Vert\Vert\int_{\Omega}\psi(\omega)a_{2}G(\omega)d\mu(\omega)\Vert\\
&\leq\Vert a_{2}^{-1}\Vert\Vert\int_{\Omega}\psi(\omega)a_{1}F(\omega)d\mu(\omega)\Vert +\Vert a_{2}^{-1}\Vert\Vert\int_{\Omega}\psi(\omega)(a_{1}F(\omega)-a_{2}G(\omega))d\mu(\omega)\Vert\\
&\leq (1+\alpha)\Vert a_{2}^{-1}\Vert\Vert\int_{\Omega}\psi(\omega)a_{1}F(\omega)d\mu(\omega)\Vert\\
&\qquad+\beta\Vert a_{2}^{-1}\Vert\Vert\int_{\Omega}\psi(\omega)a_{2}G(\omega)d\mu(\omega)\Vert +\gamma\Vert a_{2}^{-1}\Vert\Vert\psi\Vert .
\end{align*}
Since
\begin{equation*}
\Vert\int_{\Omega}\psi(\omega)a_{1}F(\omega)d\mu(\omega)\Vert=\Vert a_{1}\Vert\Vert T_{F}\psi\Vert\leq\sqrt{B}\Vert a_{1}\Vert\Vert\psi\Vert,
\end{equation*}
so
\begin{align*}
\Vert\int_{\Omega}\psi(\omega)G(\omega)d\mu(\omega)\Vert &\leq\Vert a_{2}^{-1}\Vert(1+\alpha)\sqrt{B}\Vert a_{1}\Vert\Vert\psi\Vert\\
&\qquad+\beta\Vert a_{2}^{-1}\Vert\Vert a_{2}\Vert\Vert\int_{\Omega}\psi(\omega)G(\omega)d\mu(\omega)\Vert +\gamma\Vert a_{2}^{-1}\Vert\Vert\psi\Vert,
\end{align*}
then
\begin{equation*}
\Vert\int_{\Omega}\psi(\omega)G(\omega)d\mu(\omega)\Vert\leq\dfrac{\Vert a_{2}^{-1}\Vert(\Vert a_{1}\Vert(1+\alpha)\sqrt{B}+\gamma)}{1-\beta\Vert a_{2}^{-1}\Vert\Vert a_{2}\Vert}\Vert\psi\Vert.
\end{equation*}
Also,
\begin{align*}
\Vert a_{2}\Vert\Vert\int_{\Omega}\psi(\omega)G(\omega)d\mu(\omega)\Vert &\geq\Vert\int_{\Omega}\psi(\omega)a_{2}G(\omega)d\mu(\omega)\Vert\\
&\geq\Vert\int_{\Omega}\psi(\omega)a_{1}F(\omega)d\mu(\omega)\Vert -\Vert\int_{\Omega}\psi(\omega)(a_{1}F(\omega)-a_{2}G(\omega))d\mu(\omega)\Vert\\
&\geq (1-\alpha)\Vert\int_{\Omega}\psi(\omega)a_{1}F(\omega)d\mu(\omega)\Vert\\
&\qquad -\beta\Vert\int_{\Omega}\psi(\omega)a_{2}G(\omega)d\mu(\omega)\Vert -\gamma\Vert\psi\Vert\\
&\geq (1-\alpha)\Vert a_{1}^{-1}\Vert^{-1}\Vert\int_{\Omega}\psi(\omega)F(\omega)d\mu(\omega)\Vert\\
&\qquad -\beta\Vert a_{2}\Vert\Vert\int_{\Omega}\psi(\omega)G(\omega)d\mu(\omega)\Vert -\gamma\Vert\psi\Vert,
\end{align*}
then
\begin{equation*}
(1+\beta)\Vert a_{2}\Vert\Vert\int_{\Omega}\psi(\omega)G(\omega)d\mu(\omega)\Vert\geq(1-\alpha)\Vert a_{1}^{-1}\Vert^{-1}\Vert\int_{\Omega}\psi(\omega)F(\omega)d\mu(\omega)\Vert -\gamma\Vert\psi\Vert.
\end{equation*}
Since $F$ is a Riesz-type frame, so 
\begin{equation*}
\Vert\int_{\Omega}\psi(\omega)F(\omega)d\mu(\omega)\Vert\geq M\Vert\psi\Vert.
\end{equation*}
Hence
\begin{equation*}
\Vert\int_{\Omega}\psi(\omega)G(\omega)d\mu(\omega)\Vert\geq\dfrac{(1-\alpha)M\Vert a_{1}^{-1}\Vert^{-1}-\gamma}{(1+\beta)\Vert a_{2}\Vert}\Vert\psi\Vert.
\end{equation*}
Therefore
\begin{equation*}
\dfrac{(1-\alpha)M\Vert a_{1}^{-1}\Vert^{-1}-\gamma}{(1+\beta)\Vert a_{2}\Vert}\Vert\psi\Vert\leq\Vert\int_{\Omega}\psi(\omega)G(\omega)d\mu(\omega)\Vert\leq\dfrac{\Vert a_{2}^{-1}\Vert(\Vert a_{1}\Vert(1+\alpha)\sqrt{B}+\gamma)}{1-\beta\Vert a_{2}^{-1}\Vert\Vert a_{2}\Vert}\Vert\psi\Vert.
\end{equation*}
\end{proof}
%-------------------------------------------------------------------------------------------------------------------------------------------
Duo to Theorem \ref{pert1}, if $a_{1}=a_{2}=1_{\mathcal A}$, the identity element of $\mathcal A$, then the next corollary holds.
%------------------------------------------------------------------------------------------------------------------------------------------- 
\begin{corollary}\label{pert3}
Let $F:\Omega \to U$ be a continuous frame for Hilbert $C^{\ast}$-module $U$ over a unital $C^*$-algebra $\mathcal A$ with the frame bounds $A,B>0$ and the synthesis operator $T_{F}$. Also assume that $G:\Omega \to U$ is a weakly measurable mapping.\\
$(1)\;$ If there exist constants $\alpha ,\beta, \gamma\geq 0$ such that $max\lbrace\beta , \alpha+\dfrac{\gamma}{\sqrt{A}}\rbrace<1$ and
\begin{align*}
\Vert\int_{\Omega}\psi(\omega)\langle F(\omega)-G(\omega),f\rangle d\mu(\omega)\Vert\leq & \alpha\Vert\int_{\Omega}\psi(\omega)\langle F(\omega),f\rangle d\mu(\omega)\Vert\\
& +\beta\Vert\int_{\Omega}\psi(\omega)\langle G(\omega),f\rangle d\mu(\omega)\Vert +\gamma\Vert\psi\Vert ,
\end{align*}
for all $\psi\in L^{2}(\Omega ,\mathcal A)$ and $f\in U$. Then $G$ is a continuous frame for $U$.\\
$(2)\;$ If $F$ is a Riesz-type frame for $U$ and there exist constants $\alpha ,\beta, \gamma\geq 0$ such that $max\lbrace \beta , \alpha +\frac{\gamma}{M},\alpha+\dfrac{\gamma}{\sqrt{A}}\rbrace<1$ and
\begin{equation*}
\Vert\int_{\Omega}\psi(\omega)(F(\omega)-G(\omega)) d\mu(\omega)\Vert\leq\alpha\Vert\int_{\Omega}\psi(\omega)F(\omega) d\mu(\omega)\Vert +\beta\Vert\int_{\Omega}\psi(\omega)G(\omega) d\mu(\omega)\Vert +\gamma\Vert\psi\Vert ,
\end{equation*}
for all $\psi\in L^{2}(\Omega ,\mathcal A)$. Then $G$ is a Riesz-type frame for $U$.
\end{corollary}
%-------------------------------------------------------------------------------------------------------------------------------------------
In the next corollary, we show that if $F$ is a continuous frame and $\gamma=0$ in Corollary \ref{pert3}, then $G$ is a Riesz-type frame if and only if $F$ is a Riesz-type frame. 
%-------------------------------------------------------------------------------------------------------------------------------------------
\begin{corollary}
Let $F:\Omega \to U$ be a continuous for Hilbert $C^{\ast}$-module $U$ over a unital $C^*$-algebra $\mathcal A$. Assume that $G:\Omega \to U$ is a weakly measurable mapping and there exist constants $0\leq\alpha ,\beta<1$ such that
\begin{equation*}
\Vert\int_{\Omega}\psi(\omega)(F(\omega)-G(\omega)) d\mu(\omega)\Vert\leq\alpha\Vert\int_{\Omega}\psi(\omega)F(\omega) d\mu(\omega)\Vert +\beta\Vert\int_{\Omega}\psi(\omega)G(\omega) d\mu(\omega)\Vert ,
\end{equation*}
for all $\psi\in L^{2}(\Omega ,\mathcal A)$. Then $F$ is a Riesz-type frame if and only if $G$ is a Riesz-type frame.
\end{corollary}

\begin{proof}
By the Theorem \ref{pert1}, the mapping $G$ is a continuous frame for $U$. Assume that $T_{F}$ and $T_{G}$ are the synthesis operators of $F$ and $G$, respectively. It is enough to show that $R(T^{*}_{F})=R(T^{*}_{G})$. Note that by the \cite[Lemma 1.1]{XS}, we have
\begin{equation*}
ker(T_{G})\oplus R(T^{*}_{F})=L^{2}(\Omega ,\mathcal A)=ker(T_{F})\oplus R(T^{*}_{G}).
\end{equation*}
Hence we show that $ker(T_{F})=ker(T_{G})$. For each $\varphi\in L^{2}(\Omega ,\mathcal A)$, if $\varphi\in ker(T_{G})$, then by assumption, $\Vert T_{F}\varphi\Vert\leq\beta\Vert T_{F}\varphi\Vert$ and so $T_{F}\varphi =0$. This implies that $\varphi\in ker(T_{F})$ and $ker(T_{G})\subseteq ker(T_{F})$. Similarly $ker(T_{F})\subseteq ker(T_{G})$. Therefore $ker(T_{F})=ker(T_{G})$.

\end{proof}                   
%-------------------------------------------------------------------------------------------------------------------------------------------
Now, we examine the effect of duality on the perturbation of continuous frames in Hilbert $C^{\ast}$-modules.
%-------------------------------------------------------------------------------------------------------------------------------------------
\begin{theorem}\label{pert-d}
Let $F:\Omega \to U$ be a continuous frame for Hilbert $C^{\ast}$-module $U$ over a unital $C^*$-algebra $\mathcal A$ with the frame bounds $A,B>0$ and the pre-frame operator $T_{F}$. Also assume that $G:\Omega \to U$ is a dual of $F$ with bounds $C,D>0$ and $K:\Omega \to U$ is a weakly measurable mapping which satisfies the following two conditions:\\
$(1)\;\alpha:=\int_{\Omega}\Vert F(\omega)-K(\omega)\Vert^{2} d\mu(\omega)<\infty$,\\
$(2)\;\beta:=\int_{\Omega}\Vert F(\omega)-K(\omega)\Vert\Vert G(\omega)\Vert d\mu(\omega)<1$,\\
for all $f\in U$. Then $K$ is a continuous frame for $U$.
\end{theorem}

\begin{proof}
Assume that $T_{K}:L^{2}(\Omega ,\mathcal A)\longrightarrow U$ is weakly defined by
\begin{equation*}
\langle T_{K}\varphi ,f\rangle=\int_{\Omega}\varphi(\omega)\langle K(\omega),f\rangle d\mu(\omega),
\end{equation*}
for all $\varphi\in L^{2}(\Omega ,\mathcal A)$ and $f\in U$. Then
\begin{align*}
\Vert\langle T_{K}\varphi ,f\rangle\Vert &=\Vert\int_{\Omega}\varphi(\omega)\langle K(\omega),f\rangle d\mu(\omega)\Vert\\
& \leq\Vert\int_{\Omega}\varphi(\omega)\langle F(\omega)-K(\omega),f\rangle d\mu(\omega)\Vert+\Vert\int_{\Omega}\varphi(\omega)\langle F(\omega),f\rangle d\mu(\omega)\Vert\\
& =\Vert\langle\varphi ,\lbrace\langle f,F(\omega)-K(\omega)\rangle\rbrace_{\omega\in\Omega}\rangle\Vert+\Vert\langle T_{F}\varphi ,f\rangle\Vert\\
& \leq\Vert\varphi\Vert\Vert\int_{\Omega}\langle f,F(\omega)-K(\omega)\rangle\langle F(\omega)-K(\omega),f\rangle d\mu(\omega)\Vert^{\frac{1}{2}}+\sqrt{B}\Vert\varphi\Vert\Vert f\Vert\\
&  \leq\Vert\varphi\Vert(\int_{\Omega}\Vert F(\omega)-K(\omega)\Vert^{2} d\mu(\omega))^{\frac{1}{2}}\Vert\langle f,f\rangle\Vert^{\frac{1}{2}}+\sqrt{B}\Vert\varphi\Vert\Vert f\Vert\\
& \leq(\sqrt{\alpha}+\sqrt{B})\Vert\varphi\Vert\Vert f\Vert,
\end{align*}
for all $\varphi\in L^{2}(\Omega ,\mathcal A)$ and $f\in U$. So
\begin{align*}
\Vert T_{K}\varphi\Vert & =\sup_{\Vert f\Vert\leq 1}\Vert\langle T_{K}\varphi ,f\rangle\Vert\\
& \leq\sup_{\Vert f\Vert\leq 1}(\sqrt{\alpha}+\sqrt{B})\Vert f\Vert\Vert\varphi\Vert\leq(\sqrt{\alpha}+\sqrt{B})\Vert\varphi\Vert,
\end{align*}
for all $\varphi\in L^{2}(\Omega ,\mathcal A)$ and $f\in U$. Then $T_{K}$ is well-defined and $\Vert T_{K}\Vert\leq\sqrt{\alpha}+\sqrt{B}$. This shows that $K$ is a continuous Bessel mapping for $U$.\\
Now we define $L:U\longrightarrow U$ where
\begin{equation*}
\langle Lf,g\rangle=\int_{\Omega}\langle f,G(\omega)\rangle\langle K(\omega),g\rangle d\mu(\omega),
\end{equation*}
for each $f,g\in U$. Then
\begin{align*}
\langle f-Lf,g\rangle & =\int_{\Omega}\langle f,G(\omega)\rangle\langle F(\omega)-K(\omega),g\rangle d\mu(\omega)\\
& =\langle f,\int_{\Omega}\langle g,F(\omega)-K(\omega)\rangle G(\omega) d\mu(\omega)\rangle,
\end{align*}
and
\begin{align*}
\Vert f-Lf\Vert & =\sup_{\Vert g\Vert\leq 1}\Vert\langle f-Lf ,g\rangle\Vert\\
& \leq\sup_{\Vert g\Vert\leq 1}\Vert f\Vert\Vert\int_{\Omega}\langle g,F(\omega)-K(\omega)\rangle G(\omega) d\mu(\omega)\Vert\\
& \leq\sup_{\Vert g\Vert\leq 1}\Vert f\Vert\Vert g\Vert\int_{\Omega}\Vert F(\omega)-K(\omega)\Vert\Vert G(\omega)\Vert d\mu(\omega)\leq\beta\Vert f\Vert,
\end{align*}
for all $f,g\in U$. Hence by \cite[Lemma 4.1]{RND}, $L$ is invertible and
\begin{equation*}
\Vert L\Vert\leq 1+\beta\qquad and\qquad\Vert L^{-1}\Vert\leq\dfrac{1}{1-\beta}
\end{equation*}
Then for each $f\in U$ we have
\begin{align*}
\langle f,f\rangle & =\langle LL^{-1}f,f\rangle\\
& =\int_{\Omega}\langle L^{-1}f,G(\omega)\rangle\langle K(\omega),f\rangle d\mu(\omega)\\
& =\langle\lbrace\langle L^{-1}f,G(\omega)\rangle\rbrace_{\omega\in\Omega},\lbrace\langle f,K(\omega)\rangle\rbrace_{\omega\in\Omega}\rangle,
\end{align*}
and
\begin{align*}
\Vert f\Vert^{4} &=\Vert\langle f,f\rangle\Vert^{2}\\
& \leq\Vert\lbrace\langle L^{-1}f,G(\omega)\rangle\rbrace_{\omega\in\Omega}\Vert^{2}\Vert\lbrace\langle f,K(\omega)\rangle\rbrace_{\omega\in\Omega}\Vert^{2}\\
& =\Vert\int_{\Omega}\langle L^{-1}f,G(\omega)\rangle\langle G(\omega),L^{-1}f\rangle d\mu(\omega)\Vert\Vert\int_{\Omega}\langle f,K(\omega)\rangle\langle K(\omega),f\rangle d\mu(\omega)\Vert\\
& \leq D\Vert L^{-1}f\Vert^{2}\Vert\int_{\Omega}\langle f,K(\omega)\rangle\langle K(\omega),f\rangle d\mu(\omega)\Vert\\
& \leq D\Vert L^{-1}\Vert^{2}\Vert f\Vert^{2}\Vert\int_{\Omega}\langle f,K(\omega)\rangle\langle K(\omega),f\rangle d\mu(\omega)\Vert\\
& \leq\dfrac{D}{(1-\beta)^{2}}\Vert f\Vert^{2}\Vert\int_{\Omega}\langle f,K(\omega)\rangle\langle K(\omega),f\rangle d\mu(\omega)\Vert
\end{align*}
So
\begin{equation*}
\dfrac{(1-\beta)^{2}}{D}\Vert f\Vert^{2}\leq\Vert\int_{\Omega}\langle f,G(\omega)\rangle\langle K(\omega),g\rangle d\mu(\omega)\Vert.
\end{equation*}
Therefore $K$ is a continuous frame for $U$, by Theorem \ref{c-f-norm}.
\end{proof}
%-------------------------------------------------------------------------------------------------------------------------------------------
\begin{remark}
Let $F:\Omega \to U$ be a continuous frame for Hilbert $C^{\ast}$-module $U$ over a unital $C^*$-algebra $\mathcal A$ with lower  bound $A$ and $K:\Omega \to U$ is a weakly measurable mapping such that 
\begin{equation*}
\int_{\Omega}\Vert F(\omega)-K(\omega)\Vert^{2} d\mu(\omega)<\infty.
\end{equation*}
Then $F-K$ is a continuous Bessel mapping for $U$. For,
\begin{align*}
\Vert\int_{\Omega}\langle f,& F(\omega)-K(\omega)\rangle\langle F(\omega)-K(\omega),f\rangle d\mu(\omega)\Vert\\
& \leq\int_{\Omega}\Vert\langle f,F(\omega)-K(\omega)\rangle\langle F(\omega)-K(\omega),f\rangle\Vert d\mu(\omega)\\
& \leq\int_{\Omega}\Vert F(\omega)-K(\omega)\Vert^{2} d\mu(\omega)\Vert\langle f,f\rangle\Vert ,
\end{align*}
for all $f\in U$. Assume that $N$ is the Bessel bound of $F-K$. For $K$ to be a continuous frame, condition $N<A$ must be satisfied, by Corollary \ref{pert(F-G)B}, but it is not necessary in Theorem \ref{pert-d}.
\end{remark}
%-------------------------------------------------------------------------------------------------------------------------------------------
Now, we introduce an interesting operator on Hilbert $C^{\ast}$-module $U$ and give some results about it.
%-------------------------------------------------------------------------------------------------------------------------------------------
Let $F,G:\Omega \to U$ be two continuous Bessel mappings for Hilbert $C^{\ast}$-module $U$ over a unital $C^*$-algebra $\mathcal A$ with bounds $B_{F}$, $B_{G}$ respectively. Define
\begin{align*}
R_{F,G}: &U\longrightarrow U\\
& f\longmapsto \int_{\Omega}\langle f,F(\omega)\rangle G(\omega) d\mu(\omega)
\end{align*}
for all $f\in U$. It is easily to seen that $R_{F,G}$ is well-defined, bounded, adjointable and $R_{F,G}=T_{G}T^{*}_{F}$, where $T_{F}$ and $T_{G}$ are the pre-frame operators of $F$ and $G$, respectively.
%-------------------------------------------------------------------------------------------------------------------------------------------
\begin{proposition}\label{R-S-C}
Let $F:\Omega \to U$ be a continuous frame for Hilbert $C^{\ast}$-module $U$ over a unital $C^*$-algebra $\mathcal A$ and $G:\Omega \to U$ is a continuous Bessel mapping for $U$ with the pre-frame operators $T_{F}$ and $T_{G}$, respectively. If $R_{F,G}$ is surjective, then $G$ is a continuous frame for $U$.
\end{proposition}

\begin{proof}
It is clear that $R_{F,G}=T_{G}T^{*}_{F}$. Since $R_{F,G}$ is surjective,so for each $f\in U$ there exists $h\in U$ such that $f=R_{F,G}h=T_{G}T^{*}_{F}h$. On the other hand $T^{*}_{F}h\in L^{2}(\Omega ,\mathcal A)$. Then $T_{G}$ is surjective and hence $G$ is a continuous frame for $U$, by \cite[Theorem 2.15]{GHSH}.
\end{proof}
%-------------------------------------------------------------------------------------------------------------------------------------------
\begin{remark}
The convers of the previouse proposition holds when $F$ is a Riesz-type frame. Indeed, since $F$ is a Riesz-type frame, so $T^{*}_{F}$ is onto. Also that $G$ is a continuous frame, implies that $T_{G}$ is surjective. Therefore $R_{F,G}=T_{G}T^{*}_{F}$ is surjective.
\end{remark}
%-------------------------------------------------------------------------------------------------------------------------------------------
\begin{proposition}
Let $F:\Omega \to U$ be a Riesz-type frame for Hilbert $C^{\ast}$-module $U$ over a unital $C^*$-algebra $\mathcal A$ and $G:\Omega \to U$ is a continuous Bessel mapping for $U$ with the pre-frame operators $T_{F}$ and $T_{G}$, respectively. Then $R_{F,G}$ is invertible if and only if $G$ is a Riesz-type frame for $U$.
\end{proposition}

\begin{proof}
$(\Longrightarrow)$\; 
Suppose that $R_{F,G}$ be invertible. Since $F$ is Riesz-type, so $T^{*}_{F}$ is onto and so is invertible. This shows that $T_{G}$ is also invertible. Hence $G$ is a Riesz-type frame for $U$.

$(\Longleftarrow)$\; Let $F$ and $G$ be Riesz-type frames. Then their pre-frame operators are invertible. Hence $R_{F,G}=T_{G}T^{*}_{F}$ is invertible.
\end{proof}
%-------------------------------------------------------------------------------------------------------------------------------------------
The following theorem shows if the operator frame of a continuous frame $F$ is near to $R_{F,G}$, then $G$ is a continuous frame.
%shows the effect of the approximation of the frame operator of continuous frame $F$ to the operator $R_{F,G}$ on a continuous Bessel mapping $G$.
%-------------------------------------------------------------------------------------------------------------------------------------------
\begin{theorem}\label{R-S}
Let $F:\Omega \to U$ be a continuous frame for Hilbert $C^{\ast}$-module $U$ over a unital $C^*$-algebra $\mathcal A$ with bounds $A,B$ and the continuous frame operator $S_{F}$. If $G:\Omega \to U$ is a continuous Bessel mapping for $U$ and there exists $0<\lambda<A$ such that $\Vert R_{F,G}f-S_{F}f\Vert\leq\lambda\Vert f\Vert$, then $G$ is a continuous frame for $U$.
\end{theorem}

\begin{proof}
By The proposition \ref{R-S-C}, it is enough to show that $R_{F,G}$ is surjective. For each $f\in U$ we have,
\begin{align*}
\Vert R^{*}_{F,G}f-S_{F}f\Vert &=\Vert (R_{F,G}-S_{F})^{*}f\Vert\\
& \leq\Vert (R_{F,G}-S_{F})^{*}\Vert\Vert f\Vert \\
& = \Vert R_{F,G}-S_{F}\Vert\Vert f\Vert \leq\lambda\Vert f\Vert.
\end{align*}
Also
\begin{align*}
\Vert R^{*}_{F,G}f\Vert &=\Vert R^{*}_{F,G}f-S_{F}f+S_{F}f\Vert\\
& \geq\Vert S_{F}f\Vert -\Vert R^{*}_{F,G}f-S_{F}f\Vert\geq(A-\lambda)\Vert f\Vert.
\end{align*}
Hence $R^{*}_{F,G}$ is bounded below with respect to norm and eqivalently $R_{F,G}$ is surjective, by the Lemma \ref{SB}.
\end{proof}
%-------------------------------------------------------------------------------------------------------------------------------------------
From the proof of Theorem \ref{R-S}, we know that both $R_{F,G}$ and $R^{*}_{F,G}$ are surjective, injective, closed range and so are invertible. This implies that $T_{G}$ is invertible if and only if $T^{*}_{F}$ is invertible. Hence the following corollary holds.
%-------------------------------------------------------------------------------------------------------------------------------------------
\begin{corollary}
Let $F:\Omega \to U$ be a continuous frame for Hilbert $C^{\ast}$-module $U$ over a unital $C^*$-algebra $\mathcal A$ with bounds $A,B$ and the continuous frame operator $S_{F}$. If $G:\Omega \to U$ is a continuous Bessel mapping for $U$ and there exists $0<\lambda<A$ such that $\Vert R_{F,G}f-S_{F}f\Vert\leq\lambda\Vert f\Vert$, then $F$ is a Riesz-type frame for $U$ if and only if $G$ is a Riesz-type frame for $U$.
\end{corollary}
%-------------------------------------------------------------------------------------------------------------------------------------------

\end{document}